\documentclass[12 pt]{amsart}

\usepackage{amscd, amssymb}

\newtheorem{pr}{Proposition}

\newtheorem{tm}{Theorem$^*$}

\newcommand{\proj}{\mathbf P}

\newcommand{\rarr}{\rightarrow}

\newcommand{\com}{\mathbb{C}}

\newcommand{\lan}{\langle}
\newcommand{\ran}{\rangle}
\newcommand{\eqq}{\stackrel{\sim}{=}}

\newcommand{\sss}{\Big( \frac{\sin(i\lambda/2)}{i\lambda/2} \Big)}

\def\scup{\mathbin{\text{\scriptsize$\cup$}}}
\def\scap{\mathbin{\text{\scriptsize$\cap$}}}

\newcommand{\bpf}{\noindent {\em Proof.} }

\begin{document}
\title{The Toda equations  and the Gromov-Witten theory of the
Riemann sphere}
\author{R. Pandharipande}
\date{20 December 1999}
\maketitle

\pagestyle{plain}
\setcounter{section}{-1}
\section{\bf{Introduction}}
\subsection{Toda equations}
The Gromov-Witten theory of $\proj^1$ has been intensively studied
in a sequence of remarkable papers
by Eguchi, Hori, Xiong, Yamada, and Yang  [EHY], [EY], [EYY], [EHX].
A major step in this analysis was the discovery of a 
(conjectural) matrix model for $\proj^1$. 
Let the variables $x,y$ correspond to the
unit and K\"ahler class of $\proj^1$ respectively.
Let $x_i,y_i$ (for $i\geq 0$) be the coordinates on the large phase
space. Define the potential function
\begin{equation}
\label{mmmd}
V(M) = -2M (\log(M)-1) + 2\sum_{i\geq 1}
x_i M^i \big( \log(M)- \sum_{j=1}^i \frac{1}{j} \big)
\end{equation}
$$ 
+ \sum_{i\geq 1} \frac{y_{i-1}}{i} M^i$$
on the space of $N\times N$ Hermitian matrices.
The partition function $Z$ is then defined by the
asymptotic expansions of the following Hermitian
matrix integrals:
$$Z(N, x_i, y_i) = \int_{M_{N\times N}} dM \ 
\exp(N \ tr(V(M))).$$
The Gromov-Witten potential $F=F_{\proj^1}$ is defined by
$$F(\lambda,x_i,y_i)= \sum_{g\geq 0} \lambda^{2g-2} F^g
(x_i,y_i),$$
where $F^g$ is the generating function of the 
genus $g$ descendent invariants.
Let $\gamma= \sum_{i\geq 0} x_i \tau_i(x)+y_i\tau_i(y)$. Then,
$$F^{g} = \sum_{d\geq 0} \sum_{n\geq 0} \frac{1}{n!}
\lan \gamma^n \ran_{g,d}^{\proj^1}.$$
There is a (conjectural) relationship in the large $N$
limit between
$Z$ and  $\exp (F)$
connecting  the matrix model (\ref{mmmd}) to $\proj^1$.
Related matrix models
had previously been found in the
study of two dimensional topological gravity (see [W], [K]).

A study of the matrix model for $\proj^1$
led to the discovery of the (conjectural) Toda
equations for $F$. 
The Toda equations may be written in the following
form:
\begin{equation}
\label{toda}
\exp\Big( F(x_0+\lambda) + F(x_0-\lambda)
- 2 F \Big) = \lambda^2  F_{y_0 y_0},
\end{equation}
where $F(x_0\pm \lambda)= F(\lambda, x_0\pm \lambda, x_1,\ldots,
y_0, \ldots).$
Alternatively, 
equation (\ref{toda}) is equivalent to
a sequence of differential equations indexed by genus
starting with the genus 0 and 1 equations:
\begin{equation}
\label{tzero}
\exp(F^0_{x_0 x_0}) = F^0_{y_0 y_0},
\end{equation}
\begin{equation}
\label{tone}
\exp(F^0_{x_0 x_0}) ( F^1_{x_0 x_0} + \frac{1}{12} F^0_{x_0 x_0 x_0 x_0})
= F^1_{y_0 y_0}.
\end{equation}
Properties of the Toda equations suggested the existence
of Virasoro constraints satisfied $F$. 
While the Virasoro constraints have been formulated
for arbitrary nonsingular
projective varieties [EHX] (see also [CK]), the Toda equations
seem particular to $\proj^1$.

The genus 0 equation (\ref{tzero}) is well known in several
contexts as a direct consequence of the topological
recursion relations.
A proof of the genus 1 equation is given in Section 1 --
the standard techniques of two dimensional gravity are
used. 

\begin{pr}
The Toda equations are true in genus 0 and 1.
\end{pr}

\noindent
The genus 1 result was explained to the author by E. Getzler.
The Toda equations are unproven for genus $g\geq 2$.

The main goal of this paper is to investigate the
consequences of the Toda equations. These equations
are much stronger than the Virasoro constraints for
$\proj^1$.

\begin{pr}
The Toda equations determine $F$
from the degree 0 potential $F|_{d=0}$.
\end{pr}

The degree 0 potential for $\proj^1$ has been studied via
Hodge integrals in [GeP]. In the past few years, closed form
evaluations of many of these Hodge integrals series 
have been proven
in [FaP1], [FaP2], [FaP3]. The Toda equations together
with the Hodge integral calculations lead to several explicit
predictions for basic descendent series of $\proj^1$.

\subsection{1-point invariants}
The first question considered is the computation of the
1-point invariants of $\proj^1$.
For any nonsingular projective variety $X$,
the 1-pointed genus 0 invariants play a special role in
the Gromov-Witten theory of $X$ [G] (see also [P]). A solution of the
basic differential equation in the quantum cohomology
ring $QH^*(X)$,
\begin{equation}
\label{qde}
\hbar \frac{\partial}{\partial t_i} \gamma =  \partial_i * \gamma,
\end{equation}
may be constructed from the 1-pointed genus 0 invariants --
here $\{t_i\}$ are variables indexing the divisor classes
in $H^2(X)$. 
In the case of $\proj^1$, the formulas 
$$\lan \tau_{2d-2}(y) \ran^{\proj^1}_{0, d} = \frac{1}{(d!)^2},$$
$$ \lan \tau_{2d-1}(x) \ran^{\proj^1}_{0,d} 
= - \ \frac{2 \sum_{j=1}^d \frac{1}{j}}
{(d!)^2},$$
may be derived by studying (\ref{qde}).
It is natural to study the higher genus 1-point series
for $\proj^1$.
$$Y_d(\lambda) = \sum_{g\geq 0} \lambda^{2g} \lan \tau_{2g+2d-2}(y) 
\ran^{\proj^1}_{g,d}\ ,$$
$$X_d(\lambda) = \sum_{g\geq 0} \lambda^{2g} 
\lan \tau_{2g+2d-1}(x) \ran^{\proj^1}_{g,d}.$$
In Section 2, the Toda equations and the Hodge integral evaluations
are shown to yield a conjectural computation of these series.

\begin{tm} For $d>0$, $Y_d(\lambda)$ and $X_d(\lambda)$ are determined
by:
$$ Y_d(\lambda) = \frac{1}{(d!)^2} \sss ^{2d-1},$$
$$ X_d(\lambda) =  \frac{2}{(d!)^2} \sss ^{2d-1}  \Big( \log\sss  -
\sum_{j=1}^d \frac{1}{j}\Big) .$$
\end{tm}

\noindent {\em Note:} All results in the paper marked with
an asterisk depend upon the
conjectural Toda equations.

{\vspace{+10pt}}

When $X$ is a point, the 1-point invariants have 
been computed:
$$\sum_{g\geq 0} t^{2g} \lan \tau_{3g-2} \ran_g^X = \exp(t/24),$$
with $\lan \tau_{-2} \ran _0 =1$ by definition (see [FaP1]).
Closed series forms of the 1-point invariants are not known to the
author
in any other non-trivial cases.

\subsection{Degree 1 invariants and Hodge integrals}
The degree 1 invariants $\lan \tau_{a_1}(y) \cdots \tau_{a_n}(y) 
\ran^{\proj^1}_{g,1}$
are studied via the Toda equation in Section 3. For dimension
reasons, $\sum_{i=1}^n a_i =2g$ for these degree 1 invariants.
 Let
\begin{equation}
\label{lll}
L(y_i)= \sum_{g\geq 0} \lambda^{2g-2} \sum_{n\geq 0} \frac{1}{n!}
\sum_{(a_1, \ldots, a_n), \sum a_i=2g} y_{a_1} \cdots y_{a_n} 
\lan \tau_{a_1}(y) \cdots \tau_{a_n}(y) 
\ran^{\proj^1}_{g,1}
\end{equation}
denote the generating function of the invariants.
A virtual localization argument shows:
\begin{equation}
\label{hoho}
\lan \tau_{a_1}(y) \cdots \tau_{a_n}(y) \ran_{g,1}^{\proj^1} =
\int_{\overline{M}_{g,n+1}} \prod_{i=1}^n\psi_i^{a_i}
 \cdot \frac{1-\lambda_1 + \lambda_2 - \cdots +(-1)^g \lambda_g}{1-\psi_0}.
\end{equation} 
Here, $\lambda_i$ denotes the $i^{th}$ Chern class of
the Hodge bundle on the moduli space of Deligne-Mumford
stable pointed curves $\overline{M}_{g,n+1}$.
Hence, $L(y_i)$ may be viewed as a generating function of the
Hodge integrals of type (\ref{hoho}) where $\sum_{i=1}^n{a_i}=2g$.

Let the constants $c_{2k}$ be defined by:
\begin{equation}
\label{cdeff}
\frac{\text{sin}(i\lambda/2)}{i\lambda/2} = 
\sum_{k\geq 0} c_{2k} \lambda^{2k}, \ \ c_{2k} = \frac{1}{2^{2k}(2k+1)!}.
\end{equation}

\begin{tm}
$L(y_i)$ is determined by the equation:
$$ \text{\em exp}\Big( \sum_{k\geq 0} \frac{y_{2k} \lambda^{2k}}{2^{2k}(2k+1)!}
\Big)
= \lambda^2  L(y_i).$$
\end{tm}

\noindent Or, equivalently,
$$\lan \tau_{2k_1}(y) \cdots \tau_{2k_n}(y) \ran_{g,1}^{\proj^1} =
\prod_{i=1}^n c_{2k_i}$$
when $\sum_{i} 2k_i = 2g$.

A basic vanishing result for Hodge integrals 
may be obtained as a consequence  of
Theorem$^*$ 2. 

\noindent {\bf Corollary 2.} {\em The Hodge integral
$$
\int_{\overline{M}_{g,n+1}} \prod_{i=1}^n\psi_i^{a_i}
 \cdot \frac{1-\lambda_1 + \lambda_2 - \cdots +(-1)^g \lambda_g}{1-\psi_0}$$
vanishes in the following cases
\begin{enumerate}
\item[(i)] $\sum_{i} a_i < 2g,$
\item[(ii)] $\sum_{i} a_i =2g$ and $a_j$ is odd for some j.
\end{enumerate}}
{\vspace{+10 pt}}

\noindent Theorem$^*$ 2 immediately implies (ii) of the 
Corollary. The dilaton
equation and (ii) together imply (i).
Theorem$^*$ 2 therefore yields a calculation of the Hodge integrals
(\ref{hoho}) in the first non-vanishing case.
A proof of the Corollary using localization techniques is
given in Section 3. 

\subsection{Hurwitz numbers}
The Toda equation
may be used to study the
classical (simple) Hurwitz numbers $H_{g,d}$.
$H_{g,d}$ is 
the number of nonsingular, genus $g$ 
curves expressible as  $d$-sheeted covers
of $\proj^1$ with a {\em fixed} general branch divisor in $\proj^1$.
By the Riemann-Hurwitz formula, the degree of the branch divisor
is $2g+2d-2$.
Define the generating function $H$ by:
$$H(\lambda,y_0) = \sum_{g\geq 0} \sum_{d>0} \lambda^{2g-2} e^{dy_0} 
\frac{H_{g,d}}
{(2g+2d-2)!}.$$
The Hurwitz numbers were first computed in [Hu] by 
using the combinatorics of the symmetric group. 

A result needed to connect the Hurwitz numbers
to the Toda equations is proven in Section 4.
\begin{pr}
\label{yay}
For all $g \geq 0$ and $d>0$,
$$H_{g,d} = \lan \tau_1(y)^{2g+2d-2} \ran^{\proj^1} _{g,d}.$$
\end{pr}
{\vspace{+10pt}}

\noindent
By Proposition \ref{yay}, there is an equality
$$H(\lambda, y_0)=  F|_{x_i=0, \ y_{1}=1, \ y_{i\geq 2}=0}.$$
One may therefore hope to find the Toda equations
control the Hurwitz numbers. In Section 4, the
following result is derived from the Toda equations
satisfied by $F$.

\begin{tm}
\label{qqq} $H$ satisfies the following differential equation:
\begin{equation}
\label{tto}
 \exp \Big( 
H(y_0+\lambda) + H(y_0-\lambda)-2H \Big) = \lambda^2 e^{-y_0} H_{y_0y_0}.
\end{equation}
\end{tm}

\noindent
As before, $H(y_0\pm \lambda)= H(\lambda, y_0\pm \lambda)$.
It is easily seen Theorem$^*$ 3 uniquely determines the
function $H$. Equation (\ref{tto}) resembles the
original Toda equation (\ref{toda}), but notice the shifted variable
in the argument of the exponential is now $y_0$ !

Theorem$^*$ \ref{qqq} determines the Hurwitz numbers by an elementary
recursion relation.
Let $\xi$ be an ordered sequence of vectors  
$$\xi = \Big( (g_1,d_1,k_1), (g_2,d_2,k_2), \ldots, 
(g_{l(\xi)},d_{l(\xi)},k_{l(\xi)}) \Big)$$
where $l(\xi)$ is the length of the sequence and
$(g_i\geq 0, d_i>0, k_i>0)$.
Let $P(g,d)$ be the set 
$$\{ \xi \ | \ \sum_{i} d_i =d-1,
\ \sum_{i} g_i+k_i= g+l(\xi)
\}.$$
Then, Theorem$^*$ \ref{qqq} yields:
\begin{equation}
\label{reccc}
d^2 H_{g,d} = 
\end{equation}
$$\sum_{\xi \in P(g,d)} \frac{2^{l(\xi)}}{l(\xi)!}
\binom {2g+2d-2}{2g_1+2d_1-2,2k_1, \ldots, 2g_l+2d_l-2, 2k_l}
\prod_{i=1}^{l(\xi)} d_i^{2k_i} H_{g_i,d_i}.$$
This conjectured form is  certainly the simplest recursion for the Hurwitz
numbers known to the author. It would be very interesting to find
a geometric or combinatorial proof of (\ref{reccc}).

We note 
the Hurwitz numbers are also related to Hodge integrals.
In [FanP],
the following formula is proven:
\begin{equation}
\label{hodgehur}
H_{g,d}= \frac{(2g+2d-2)!}{d!} \int_{\overline{M}_{g,d}}
\frac{1-\lambda_1+\lambda_2 -\lambda_3 + \ldots + (-1)^g \lambda_g}
{\prod_{i=1}^d (1-\psi_i)}
\end{equation}
(see also the announcement [ELSV]). Generating functions
for these Hurwitz numbers have been studied in [GJV] via
the relationship to Hodge integrals. However, Theorem$^*$ 3 appears
quite different from the fixed genus recursions found in
[FanP], [GJV].

\subsection{Acknowledgments}
Many thanks are due to C.-S. Xiong for explaining the
Toda equation to the author 
during a visit to the International Center for Theoretical
Physics in Trieste in the summer of 1999. 
The proof of the Toda equation in genus 1 was shown to
the author by E. Getzler. The author has also
benefitted from conversations with J. Kock and R. Vakil. 
This research was partially supported by
DMS-9801574 and an A. P. Sloan foundation
fellowship.

\section{The Toda equations in genus 0 and 1}
\subsection{Genus 0}
The Toda equation in genus 0 is:
\begin{equation}
\label{rrrr}
\exp(F^0_{x_0 x_0}) = F^0_{y_0 y_0}.
\end{equation}
The proof of (\ref{rrrr}) is very well known -- first
given in [DW].
The restriction of $F^0$ to the small phase space is:
$$F^0|_{x_{i\geq1}=0, \ y_{i\geq 1}=0} = \frac{1}{2} x_0^2y_0 + e^{y_0}.$$
Hence,
equation (\ref{rrrr}) is true when restricted to the small phase space.
The topological recursion relations in genus 0 imply:
$$F^0_{y_0y_0t_k} = F^0_{y_0y_0\bullet} g^{\bullet *} F^0_{*t_{k-1}},$$
$$\frac{\partial}{\partial t_k} \exp(F^0_{x_0x_0}) =
\frac{\partial}{\partial \bullet} \exp(F^0_{x_0x_0}) g^{\bullet *} 
F^0_{*t_{k-1}},$$
where $t$ equals $x$ or $y$, and $$\bullet g^{\bullet*} *$$
is the standard diagonal splitting in $\proj^1\times \proj^1$.
As these two recursion relations have the same form, we
conclude equation (\ref{rrrr}) is true on the large phase space.
The genus 0 Toda equation is therefore proven.

\subsection{Genus 1}
The Toda equation in genus 1 is:
\begin{equation}
\label{wwww}
\exp(F^0_{x_0x_0}) ( F^1_{x_0 x_0} + \frac{1}{12} F^0_{x_0 x_0 x_0 x_0})
= F^1_{y_0 y_0}.
\end{equation}
The proof of (\ref{wwww}) was explained to the author by E. Getzler.

Start by defining the following functions:
\begin{equation*}
A_0=F^0_{x_0x_0}, \ \ A_i= \frac{\partial ^i A_0}{\partial x_0^i},
\end{equation*}
$$B_0=F^0_{x_0y_0}, \ \ B_i= \frac{\partial ^i B_0}{\partial x_0^i}.$$
Let $Q= \exp(F_{x_0x_0})$.
Let ${\mathbb{A}}=\com[A_i,B_i,Q]$ be the 
polynomial ring  in the variables $\{A_i,B_i,Q\}$ --
no algebraic relations are imposed.
Via the genus 0 Toda equation,
differentiation by $x_0$ and $y_0$ define
unique linear transformations
\begin{equation}
\label{nnn} 
\frac{\partial}{\partial x_0}, \ \frac{\partial}{\partial y_0} 
\ : {\mathbb{A}} \rarr {\mathbb{A}}.
\end{equation}
For example,
$$\frac{\partial B_0}{\partial {y_0}}= \frac{\partial Q}{\partial x_0} 
= Q A_1.$$
Next, define the function  $\bigtriangleup$ by:
$$\bigtriangleup =  (F^0_{x_0x_0y_0})^2 - F^0_{x_0x_0x_0}F^0_{x_0y_0y_0}
= B_1^2 - Q A_1^2.$$
The differential operations (\ref{nnn})
extend uniquely to the localized ring ${\mathbb{A}}_\bigtriangleup$ (in
which $\bigtriangleup^{-1}$ is adjoined).
We will reduce (\ref{wwww}) to a computation in 
${\mathbb{A}}_\bigtriangleup$.

The following expression for the genus 1 potential of $\proj^1$
is proven in [DW] (see also [DZ], [Ge]):
$$F^1 = -\frac{1}{24} F^0_{x_0x_0} + \frac{1}{24} \log ( \bigtriangleup).$$
Therefore, we find
$$(F^1 + \frac{1}{12} F^0_{x_0x_0})_{x_0x_0} = \frac{1}{24}( A_0 +
\log \bigtriangleup) _{x_0x_0}.$$
On the other hand,
$$F^1_{y_0y_0} = 
\frac{1}{24}(-A_0 + \log ( \bigtriangleup)) _{y_0y_0}.$$
Finally, the Toda equation in genus 1 is equivalent to:
\begin{equation}
\label{denn}
Q (A_0 + \log \bigtriangleup)_{x_0x_0} = (-A_0
+ \log \bigtriangleup)_{y_0y_0}.
\end{equation}
The above  equation can be expanded  completely in 
${\mathbb{A}}_\bigtriangleup$.
A direct algebraic check shows the two sides of (\ref{denn})
are equal in ${\mathbb{A}}_\bigtriangleup$. 
This finishes the proof of the Toda equation in genus 1
and the proof of Proposition 1.

\subsection{Proposition 2}
The Toda equation may be written in the
following form:
\begin{equation}
\label{wer}
\exp \Big( \sum_{k>0} \frac{2 \lambda^{2k}}{(2k)!} \frac{\partial^{2k} F}
{\partial x_0^{2k}} \Big) = \lambda^2 F_{y_0y_0}.
\end{equation}
An important property of the argument of the exponential is:
\begin{equation}
\label{ttyy}
\sum_{k>0} \frac{2 \lambda^{2k}}{(2k)!} \frac{\partial^{2k} F}
{\partial x_0^{2k}} |_{x_i=0, \ y_{i\geq 1}=0} \ = \ y_0.
\end{equation}
Let $$D= \frac{\partial^{n+m}}
{ \partial x_{a_1}\cdots \partial x_{a_n} \partial
y_{b_1} \cdots \partial y_{b_m}}$$ 
be a differential operator.
By the divisor equation,
\begin{equation}
\label{qwett}
\lambda^2 D F_{y_0y_0}|_{x_i=0, \ y_{i\geq 1}=0} \ = 
\end{equation}
$$
\sum_{g\geq 0} \lambda^{2g} \sum_{d\geq 0}
e^{dy_0}\Big(d^2 \lan\prod_{i} \tau_{a_i}(x) \prod_j \tau_{b_j}(y) \ran
_{g,d}^{\proj^1} + \text{ lower terms}\Big).$$
The {\em lower terms} in the above equation refer to
degree $d$ 
descendent corrections via the divisor equation with {\em strictly fewer} 
$\tau_*(x)$ factors. In particular, if $n=0$, there
are no correction terms. By the Toda equation and (\ref{ttyy}),
the series (\ref{qwett}) equals
\begin{equation}
\label{lhl}
\tilde{D}\Big( 
\sum_{k>0} \frac{2 \lambda^{2k}}{(2k)!} \frac{\partial^{2k} F}
{\partial x_0^{2k}} \Big ) \cdot e^{y_0} |_{x_i=0, \ y_{i\geq 1}=0},
\end{equation}
where $\tilde{D}$ is differential expression 
determined by $D$. The factor of $e^{y_o}$ in (\ref{lhl})
implies degree $d-1$ invariants in (\ref{lhl}) are matched with
degree $d$ invariants in (\ref{qwett}).
An induction on the number of $\tau_*(x)$
factors then proves the Toda equation determines $F$
from the degree 0 potential $F|_{d=0}$.
Proposition 2 is proven.
\subsection{Degree 0}
There are two series of non-vanishing 
degree 0 descendent invariants of
$\proj^1$ [GeP]:
$$\lan \tau_{a_1}(x) \ldots \tau_{a_n}(x) \tau_b(y) \ran_{g,0}^{\proj^1}\ , 
\ \
\lan \tau_{a_1}(x) \ldots, \tau_{a_n}(x) \ran_{g,0}^{\proj^1}.$$
Both the above series are uniquely determined by the (conjectural) Virasoro
constraints from the basic  1-point series calculated via Hodge integrals
in [FaP1]:
\begin{equation}
\label{psss}
1+ \sum_{g\geq 1} \lambda^{2g} \lan \tau_{2g-2}(y) 
\ran_{g,0}^{\proj^1} = \sss^{-1},
\end{equation}
see [GeP].
The Virasoro prediction [GeP]:
$$\lan \tau_{a_1}(x) \ldots \tau_{a_n}(x) \tau_b(y) \ran_{g,0}^{\proj^1} =
\binom{2g-2+n}{a_1, \ldots,a_n,b} \cdot \lan \tau_{2g-2}(y)\ran _{g,0}
^{\proj^1}$$
has been proven in [FaP2]. While
the Virasoro prediction for the second degree 0 series
is an explicit recursion, no simple closed form
has yet been found.
However, the $1$-point
formula
\begin{equation}
\label{tsss}
\sum_{g\geq 1} \lambda^{2g} \lan  \tau_{2g-1}(x)
\ran_{g,0}^{\proj^1} = 2 \sss^{-1} \log \sss .
\end{equation}
has been proven in [FaP1].
These results and Virasoro predictions
together give good control of the degree 0 descendent invariants
of $\proj^1$. 
In contrast, the Toda equations
yield no information about the degree 0 invariants.

\section{1-point invariants of $\proj^1$}
The $1$-point series of $\proj^1$ 
will be used in the following forms:
\begin{equation}
\label{pss}
Y_0=\sum_{g\geq 0} \lambda^{2g} \lan \tau_0(x)^2\tau_{2g}(y) 
\ran_{g,0}^{\proj^1} = \sss^{-1},
\end{equation}
\begin{equation}
\label{tss}
X_0=\sum_{g\geq 0} \lambda^{2g} \lan \tau_0(x)^2 \tau_{2g+1}(x)
\ran_{g,0}^{\proj^1} = 2 \sss^{-1} \log \sss .
\end{equation}
The $\tau_0(x)^2$ terms are included to for proper treatment
of the $g=0$ term. Via the string equation, these results
are identical to (\ref{psss} -\ref{tsss}).

A simple extraction of the terms $y_{2g+2d-2}e^{dy_0}$ (corresponding to
1-point invariants) in the
Toda equation in form (\ref{wer}) yields the recursion: $d>0$,
$$\Big( \sum_{k>0}\frac{2 \lambda^{2k-2} }{(2k)!}
\Big) 
Y_{d-1} = d^2 Y_d.$$
There is an equality
$$\sum_{k>0}\frac{2 \lambda^{2k-2} }{(2k)!} = \sss ^2.$$
Therefore, the Toda equation together with (\ref{pss})
yields:
\begin{equation}
\label{dfg}
Y_d = \frac{1}{(d!)^2} \sss^{2d-1}.
\end{equation}

Similarly,
extraction of the terms $x_{2g+2d-1}e^{dy_0}$ together
with the divisor equation in Gromov-Witten theory yield:
$d>0$,
$$ \sss ^2 
X_{d-1} = d^2 X_d + 2d Y_d.$$
Using (\ref{dfg}), we find
$$ \frac{1}{d^2} 
\sss^2 X_{d-1} -  \frac{2}{d} \frac{1}{(d!)^2} \sss^{2d-1}  = X_{d}.$$
The closed form 
$$ X_d(\lambda) = \frac{2}{(d!)^2} \sss ^{2d-1}  \Big( \log\sss  - 
\sum_{j=1}^d \frac{1}{j}\Big) .$$
then follows by induction starting from (\ref{tss}).
The derivation of Theorem$^*$ 1 is complete.

\section{Degree 1 invariants of $\proj^1$ and  Hodge integrals}
\subsection{Theorem$^*$ 2}
The generating function (\ref{lll}) of the
invariants 
$$\lan \tau_{a_1}(y) \ldots \tau_{a_n}(y) \ran_{g,1}^{\proj^1}$$
is related to the full potential function by:
$$L(y_i) = e^{y_0} \cdot \text{Coeff}(e^{y_0}, \ F|_{x_{i\geq 0}=0}).$$
By the Toda equation,
$$\lambda^2 L(y_i)= e^{y_0} \cdot \text{Coeff} \Big( e^{y_0}, \ \exp
\Big( \sum_{k>0} \frac{2 \lambda^{2k}}{(2k)!} \frac{\partial^{2k} F}
{\partial x_0^{2k}} \Big)|_{x_{i\geq 0} =0} \ \Big),$$
$$= \exp \Big( \sum_{g\geq 0, k>0} 
\frac{2 \lambda^{2k-2}}{(2k)!} \ \lambda^{2g} y_{2g+2k-2} \ \lan 
\tau(x)^2_0\tau_{2g}(y) \ran_{g,0}^{\proj^1} \Big).$$
The string and divisor equations were used in the 
above equalities.
Recall the series evaluations:
$$\sum_{k>0} \frac{2 \lambda^{2k-2}}{(2k)!} = \sss^2,$$
$$ \sum_{g\geq 0} \lambda^{2g} \lan \tau_{0}(x)^2 \tau_{2g}(y)\ran
_{g,0}^{\proj^1} = \sss^{-1}.$$
Since
$$\sum_{k>0} \frac{2 \lambda^{2k-2}}{(2k)!} \cdot 
\sum_{g\geq 0} \lambda^{2g} \lan \tau_{0}(x)^2 \tau_{2g}(y)\ran
_{g,0}^{\proj^1} = \sss,$$
we conclude:
$$\lambda^2 L(y_i) = \exp ( \sum_{k\geq 0} c_{2k} y_{2k} \lambda^{2k})$$
using the definition (\ref{cdeff}) of the constants $c_{2k}$.
The derivation of Theorem$^*$ 2 is complete. 

\subsection{Hodge integrals}
\label{hhd}
The Hodge integral expression,
\begin{equation}
\label{hohoho}
\lan \tau_{a_1}(y) \cdots \tau_{a_n}(y) \ran_{g,1}^{\proj^1} =
\int_{\overline{M}_{g,n+1}} \prod_{i=1}^n\psi_i^{a_i}
 \cdot \frac{1-\lambda_1 + \lambda_2 - \cdots +(-1)^g \lambda_g}{1-\psi_0},
\end{equation} 
is proven by an application of the virtual localization formula
[GrP] to the Gromov-Witten invariants of $\proj^1$.
Equip $\proj^1$ with the standard linear torus
action with 2 fixed points $p_1, p_2$.
If the K\"ahler class is linearized to have weight 0 over
$p_1$, then the markings on all fixed point loci
contributing to (\ref{hohoho}) in the localization graph sum
must lie over $p_2$. 
The vanishing
$$\int_{\overline{M}_{h,1}} \frac{1-\lambda_1 + \lambda_2 - 
\cdots +(-1)^h \lambda_h}{1-\psi_0} =0$$
for all $h>0$ is proven in [FaP1] (and also follows from
(\ref{hodgehur})).
This vanishing implies graphs contributing to 
(\ref{hohoho}) may not have marking-free vertices of positive
genus lying over a fixed point.
Taken together, the marking and genus restrictions show
that only one  graph  has a nonzero contribution
to (\ref{hohoho}). The contributing graph is the unique
degree 1 graph with one vertex of genus $g$
carrying the all markings and lying over $p_2$
(see [GrP], [FaP1]).
The Hodge integral on the right of (\ref{hohoho})
is the
vertex integral obtained from the unique graph.

\subsection{Corollary 2}. Hodge integral relations
via virtual localization 
(as pursued in [FaP1], [FaP2], [FaP3]) easily imply
the vanishings of Corollary 2.

First, the vanishing (ii) is proven. Let the genus $g$ be fixed and
let $(a_1, \ldots, a_n)$ be a sequence of non-negative
integers with sum $|a|<2g$.
Then, the
integral
\begin{equation}
\label{dfdff}
\int_{[\overline{M}_{g,n}(\proj^1, 1)]^{vir}}  \text{ev}_1^*(y^{2g-|a|})
\scup
\prod_{i=1}^n \psi_i^{a_i} \scup
\text{ev}_i^*(y) 
\end{equation}
is well-defined and clearly vanishes as $\text{ev}_1^*(y^2)=0$.
However, when computed by virtual localization with the
linearization specified in Section \ref{hhd}, we find
(\ref{dfdff}) equals 
$$
\int_{\overline{M}_{g,n+1}} \prod_{i=1}^n\psi_i^{a_i}
 \cdot \frac{1-\lambda_1 + \lambda_2 - \cdots +(-1)^g \lambda_g}{1-\psi_0}.
$$
The proof of (ii) is complete.

To prove (i), let $(a_1, \ldots, a_n)$ be a sequence with
sum $|a|=2g$ and $a_j$ odd. We will compute the descendent integral
(\ref{dfdff}) via localization as specified in Section \ref{hhd}
with one change. Linearize the factor
$\text{ev}^*_j(y)$ to have  weight 0
over $p_2$. Then, in the localization graph sum computing
(\ref{dfdff}), it is easily seen that 
the vanishing (ii) may be applied to one vertex in every graph. Hence,
the vanishing (i) is proven and the proof of Corollary 2 is
complete.

\section{Hurwitz numbers}
\subsection{Proposition 3}
Our first goal is to prove Proposition 3 relating 
the Hurwitz numbers to descendent invariants. 
For all $g \geq 0$ and $d>0$,
$$H_{g,d} = \lan \tau_1(y)^{2g+2d-2} \ran^{\proj^1} _{g,d}.$$

\bpf Let $n=2g+2d-2$. By the hypotheses, $n\geq 0$ (and
$n=0$ only in the trivial case $g=0, d=1$).
Let $p_1, \ldots, p_n$ be fixed general points in $\proj^1$.
Consider the moduli space of stable maps 
$\overline{M}_{g,n}(\proj^1,d)$
and the closed subspace
$$ V = \text{ev}_1^{-1}(p_1) \scap \cdots \scap \text{ev}_n^{-1}(p_n)
\subset \overline{M}_{g,n}(\proj^1,d),$$
where $\text{ev}_i$ denotes the $i^{th}$ evaluation map.

Let $L_i$ denote the restriction of the $i^{th}$ cotangent
line to $V$. For each $i$, there is a canonical section 
$$s_i \in H^0(V, L_i)$$
(defined up to
a scalar $\com^*$). The section $s_i$ is obtained by the following
construction.
Let $\com \eqq T^*_{p_i}(\proj^1)$.
Let $[\mu:C \rarr \proj^1, m_1, \ldots, m_n] \in V$.
There is a canonical differential map
\begin{equation}
\label{jjj}
d\mu^*: \com \eqq T^*_{p_i}(\proj^1) \rarr L_i|_{[\mu]}.
\end{equation}
The map (\ref{jjj}) determines the section $s_i$.

An easy inductive argument proves the following set theoretic
claim:
the common zero locus $Z\subset V$ of the sections $s_1, \ldots, s_n$
equals the set of Hurwitz covers of nonsingular,
irreducible curves simply ramified over
$p_1, \ldots, p_n$. 
If $[\mu] \in Z$ has corresponds to a map with
a nonsingular domain curve $C$, then
$\mu$ must be a simply ramified Hurwitz cover by the
Riemann-Hurwitz formula.
If $[\mu] \in Z$ is nodal or irreducible, then a contradiction
is reached by counting the nodes and applying the
inductive assumption to the components of
the normalization of the domain curve.

The descendent  $\lan \tau_1(y)^{2g+2d-2} \ran_{g,d}^{\proj^1}$
may be computed as
an intersection against the virtual class of the moduli space
of maps:
\begin{equation}
\label{yyyyy}
\lan \tau_1(y)^{2g+2d-2} \ran_{g,d}^{\proj^1} =
[V]^{vir} \cap \prod_{i=1}^n (s_i),
\end{equation}
where $(s_i)$ denotes the zero locus of $s_i$.
In [FanP], it is proven the open moduli space 
$$M_{g,n}(\proj^1,d) \subset \overline{M}_{g,n}(\proj^1,d)$$
is a nonsingular Deligne-Mumford stack of expected dimension.
By Bertini's theorem, $V \scap M_{g,n}(\proj^1,d)$ is
also a nonsingular stack of expected dimension.
Hence, in an open set containing the intersection $Z$,
$[V]^{vir}$ is simply the ordinary fundamental class of $V$.

The required
transversality of the intersection cycle $Z$
follows easily from the classical fact that the space of pointed
Hurwitz covers is \'etale over the space of
branch points. We have proven
(\ref{yyyyy}) is a transverse intersection in the smooth
Deligne-Mumford
stack $V$
exactly counting
the Hurwitz covers (with the natural automorphism factors
from the orbifold geometry).
Proposition 3 is proven.

In fact, results relating first descendents to enumerative
tangency
conditions can be pursued more generally (this topic
will be discussed in an upcoming paper with T. Graber and
J. Kock). See [V] for a result closely related to Proposition 3.

\subsection{Theorem$^*$ 3}
The generating function  $H(\lambda, y_0)$ for the Hurwitz numbers
is related to the full potential function by:
$$ H(\lambda,y_0) = F|_{x_i=0,\ y_1=1, \ y_{i\geq 2}=0}.$$
By the Toda equation,
$$\lambda^2 H_{y_0y_0}= \exp 
\Big( \sum_{k>0} \frac{2 \lambda^{2k}}{(2k)!} \frac{\partial^{2k} F}
{\partial x_0^{2k}}|_{x_i=0,\ y_1=1, \ y_{i\geq 2}=0} \Big).$$
Consider first the $k=1$ term of the argument of the exponential:
\begin{equation}
\label{qqwq}
\sum_{g\geq 0 } \frac{2 \lambda^{2}}{2!} \lambda^{2g-2}  F^g_{x_0x_0}
|_{x_i=0,\ y_1=1, \ y_{i\geq 2}=0} =
\end{equation}
$$
y_0 + \sum_{g\geq 0}\sum_{d>0} \frac{2 (d\lambda)^2}{2!} \lambda^{2g-2} 
e^{dy_0} 
\frac{H_{g,d}} {(2g+2d-2)!}.
$$
The combinatorics of the
string and divisor equations are used in (\ref{qqwq}) --
together with Proposition 3.
There is a special $d=0$ term $y_0$ in (\ref{qqwq}).
The $k>1$ terms of the argument have a uniform expansion:

\begin{equation}
\label{dffd}
\sum_{g\geq 0}
\frac{2 \lambda^{2k}}{(2k)!} \lambda^{2g-2} \frac{\partial^{2k} F^g}
{\partial x_0^{2k}}|_{x_i=0,\ y_1=1, \ y_{i\geq 2}=0} =
\end{equation}
$$\sum_{g \geq 0}\sum_{d>0} \frac{2 (d\lambda)^{2k}}{(2k)!} \lambda^{2g-2}
e^{dy_0} \frac{H_{g,d}} {(2g+2d-2)!}.$$
Equations  (\ref{qqwq}) and (\ref{dffd})
together yield:
\begin{equation}
\label{fred}
\lambda^2 H_{y_0y_0}= e^{y_0} \exp 
\Big( \sum_{g\geq 0} \sum_{d>0} \sum_{k>0} 
\frac{2 (d\lambda)^{2k}}{(2k)!} 
\lambda^{2g-2}
e^{dy_0} \frac{H_{g,d}} {(2g+2d-2)!}
\Big).
\end{equation}

There is a  series equality
$$\sum_{k>0} \frac{2 (d\lambda)^{2k}}{(2k)!} =
 e^{d\lambda} + e^{-d\lambda} - 2.$$
Substituting this series into (\ref{fred}) yields:
$$\lambda^2 H_{y_0y_0}= e^{y_0} \exp 
\Big( \sum_{g\geq 0} \sum_{d>0} 
(e^{d\lambda} + e^{-d\lambda} - 2)
\lambda^{2g-2}
e^{dy_0} \frac{H_{g,d}} {(2g+2d-2)!}
\Big).$$
The final form,
$$\lambda^2 H_{y_0y_0}=  e^{y_0} \exp\Big( 
H(y_0+\lambda) + H(y_0-\lambda)-2H \Big),$$
completes the derivation of Theorem$^*$ 3.
The explicit recursion (\ref{reccc}) follows directly
from equation (\ref{fred}).

\vspace{+10 pt}
\noindent
Department of Mathematics \\
\noindent California Institute of Technology \\
\noindent Pasadena, CA 91125 \\
\noindent rahulp@cco.caltech.edu
\end{document}